# LAWS OF LARGE NUMBERS FOR EPIDEMIC MODELS WITH COUNTABLY MANY TYPES

By A. D. Barbour[1] and M. J. Luczak[2]

*Universität Zürich and London School of Economics*

In modeling parasitic diseases, it is natural to distinguish hosts according to the number of parasites that they carry, leading to a countably infinite type space. Proving the analogue of the deterministic equations, used in models with finitely many types as a "law of large numbers" approximation to the underlying stochastic model, has previously either been done case by case, using some special structure, or else not attempted. In this paper we prove a general theorem of this sort, and complement it with a rate of convergence in the $\ell_1$-norm.

**1. Introduction.** This paper is concerned with generalizations of the stochastic models introduced in Barbour and Kafetzaki (1993) and developed in Luchsinger (2001a, 2001b), which describe the spread of a parasitic disease. With such diseases, it is natural to distinguish hosts according to the number of parasites that they carry. Since it is not usually possible to prescribe a fixed upper limit for the parasite load, this leads to models with countably infinitely many types, one for each possible number of parasites. The model considered by Kretzschmar (1993) is also of this kind, though framed in a deterministic, rather than a stochastic form. Then there are models arising in cell biology, in which, for instance, hosts may be replaced by cells which are distinguished according to the number of copies of a particular gene that they carry, a number which is again, in principle, unlimited; see Kimmel and Axelrod (2002), Chapter 7, for a selection of branching process examples. The metapopulation model of Arrigoni (2003) also allows for

Received November 2006; revised January 2008.

[1]Supported in part by Schweizerischer Nationalfonds Projekte Nr. 20–107935/1 and 20–117625/1.

[2]Supported in part by the Nuffield Foundation.

*AMS 2000 subject classifications.* 92D30, 60J27, 60B12.

*Key words and phrases.* Epidemic models, infinitely many types, quantitative law of large numbers.

infinitely many types of patches, here distinguished by the size of the population in the patch.

The fact that there are infinitely many types can cause difficulty in problems which, for finitely many types, would be quite standard. To take a well known example, in a super-critical Galton–Watson branching process with finitely many types, whose mean matrix is irreducible and aperiodic and whose offspring distributions have finite variance, the proportions of individuals of the different types converge to fixed values if the process grows to infinity: a rather stronger result is to be found in Kesten and Stigum (1966). If there are infinitely many types, little is generally known about the asymptotics of the proportions, except when the mean matrix is $r$-positive, a condition which is automatically satisfied in the finite case; here, Moy (1967) was able to prove convergence under a finite variance condition. For epidemic models analogous to those above, but with only finitely many types, there are typically "law of large numbers" approximations, which hold in the limit of large populations, and are expressed in the form of systems of differential equations: see, for example, Bailey (1968) or Kurtz (1980). Proving such limits for models with infinite numbers of types is much more delicate. Kretzschmar (1993) begins with the system of differential equations as the model, and so does not consider the question; in Barbour and Kafetzaki (1993) and Luchsinger (2001a, 2001b), the arguments are involved, and make use of special assumptions about the detailed form of the transition rates.

In this paper we prove a law of large numbers approximation with explicit error rates in some generality. The models that we allow are constructed by superimposing state-dependent transitions upon a process with otherwise independent and well-behaved dynamics within the individuals; the state-dependent components are required to satisfy certain Lipschitz and growth conditions, to ensure that the perturbation of the underlying semi-group governing the independent dynamics is not too severe. The main approximation is stated in Theorem 3.1, and bounds the difference between the normalized process $N^{-1}X_N$ and a deterministic trajectory $x$ with respect to the $\ell_1$-norm, uniformly on finite time intervals. Here, $N$ is a "typical" population size, and $X_N$ is shorthand for $\{(X_N^j(t), j \in \mathbb{Z}_+), t \geq 0\}$, where $X_N^j(t)$ denotes the number of hosts at time $t$ having $j$ parasites. The theorem is sufficiently general to cover all the epidemic models mentioned above, except for that of Kretzschmar (1993), where only a modified version can be treated.

The processes that we consider can be viewed in different ways. One is to interpret them as mean-field interacting particle systems. Léonard (1990) has used this approach to study the large population behavior of a number of epidemic models, though he requires the generators to be bounded, which is an unnatural assumption in our parasite systems. Alternatively, the value $N^{-1}X_N(t)$ of the normalized process at time $t$ can be seen as a measure on $\mathbb{Z}_+$, establishing connections with the theory of measure-valued processes. Eibeck and Wagner (2003) prove law of large numbers limit theorems in a rather general setting, aimed principally at coagulation and fragmentation models, and allowing infinite numbers of types and unbounded generators; however, they do not establish error bounds.

Our argument proceeds by way of an intermediate approximation, based on a system $\widetilde{X}_N$ consisting of independent particles, which has dynamics reflecting the average behavior of $X_N$. The deterministic trajectory $x$ is discussed in Section 3, the approximation of $N^{-1}\widetilde{X}_N$ by $x$ in Section 4, culminating in Theorem 4.5, and the final approximation of $N^{-1}X_N$ by $x$ in Section 5.

**2. Specifying the model.** Our model is expressed in terms of a sequence of processes $X_N$ having state space $\mathcal{X} := \{\xi \in \mathbb{Z}_+^\infty : \sum_{i \geq 0} \xi^i < \infty\}$, where $\mathbb{Z}_+$ denotes the nonnegative integers and $\xi^i$ the $i$th component of $\xi$; we also use $e(i)$ to denote the $i$th coordinate vector. The process $X_N(t) := (X_N^j(t) : j \in \mathbb{Z}_+)$, $t \geq 0$, has $\sum_{j \geq 0} X_N^j(0) = N$, and evolves as a pure jump Markov process with transition rates given by

$$\xi \to \xi + e(j) - e(i) \text{ at rate } \xi^i\{\bar{\alpha}_{ij} + \alpha_{ij}(N^{-1}\xi)\}, \qquad i \geq 0,\ j \geq 0,\ j \neq i;$$
$$\xi \to \xi + e(i) \text{ at rate } N\beta_i(N^{-1}\xi), \qquad i \geq 0;$$
$$\xi \to \xi - e(i) \text{ at rate } \xi^i\{\bar{\delta}_i + \delta_i(N^{-1}\xi)\}, \qquad i \geq 0,$$

for any $\xi \in \mathcal{X}$, where the *nonnegative* quantities $\bar{\alpha}_{ij}, \alpha_{ij}, \beta_i, \bar{\delta}_i$ and $\delta_i$ are used to model different aspects of the underlying parasite life cycle.

We interpret $X_N^i(t)$ as the number of hosts who carry $i$ parasites at time $t$. If only the constant terms in the transitions were present, parasite communities would develop independently within different hosts, according to a pure jump Markov process with infinitesimal matrix $\bar{\alpha}$, also including the possibility of host death at rate $\bar{\delta}_i$ for hosts with $i$ parasites.

The remaining terms in the transitions are used to model interactions, and hence vary as a function of the levels $x = N^{-1}\xi \in N^{-1}\mathcal{X}$ of infection in the host population. The $\alpha_{ij}$ are associated with transitions in which the number of parasites within a particular host changes, and can be used to model interactions involving both infection with new parasites and loss of infection through parasite death. The remaining transitions allow one to model state varying rates of births, deaths and immigration of hosts, in the latter case possibly themselves infective.

The components of the transitions rates are required to satisfy a number of conditions. First, we address the $\bar{\alpha}_{ij}$ and $\bar{\delta}_i$. Letting $\Delta$ denote an absorbing "cemetery" state, reached if a host dies, set

(2.1) $\qquad \bar{\alpha}_{i\Delta} := \bar{\delta}_i, \qquad \bar{\alpha}_{ii} := -\alpha^*(i) - \bar{\delta}_i, \qquad i \geq 0,$

where we assume that $\alpha^*(i) := \sum_{j \geq 0, j \neq i} \bar{\alpha}_{ij} < \infty$. Then $\bar{\alpha}$ is the infinitesimal matrix of a time homogeneous pure jump Markov process $W$ on $\mathbb{Z}_+ \cup \Delta$. Writing

$$(2.2) \qquad p_{ij}(t) := \mathbf{P}[W(t) = j \mid W(0) = i],$$

for $i \geq 0$ and $j \in \mathbb{Z} \cup \Delta$, we shall assume that $\bar{\alpha}$ is such that $W$ is nonexplosive and that

$$(2.3) \qquad \mathbf{E}_i^0\{(W(t)+1)\} = \sum_{j \geq 0}(j+1)p_{ij}(t) \leq (i+1)e^{wt}, \qquad i \geq 0,$$

for some $w \geq 0$, where we use the notation

$$\mathbf{E}_i^0(f(W(t))) := \mathbf{E}\{f(W(t))I[W(t) \notin \Delta] \mid W(0) = i\}.$$

We shall further require that, for some $1 \leq m_1, m_2 < \infty$,

$$(2.4) \qquad \alpha^*(i) + \bar{\delta}_i \leq m_1(i+1)^{m_2} \qquad \text{for all } i \geq 0,$$

and also that, for each $j \geq 0$,

$$(2.5) \qquad \limsup_{l \to \infty} \bar{\alpha}_{lj} < \infty.$$

The remaining elements depend on the state of the system through the argument $x := N^{-1}\xi$. In the random model, $x \in N^{-1}\mathcal{X}$ has only finitely many nonzero elements, but when passing to a law of large numbers approximation, this need not be appropriate in the limit. We shall instead work within the larger spaces

$$(2.6) \qquad \ell_{11} := \left\{ x \in \mathbb{R}^\infty : \sum_{i \geq 0}(i+1)|x^i| < \infty \right\},$$

with norm $\|x\|_{11} := \sum_{i \geq 0}(i+1)|x^i|$, and the usual $\ell_1$ with norm $\|x\|_1 := \sum_{i \geq 0}|x^i|$. We then assume that $\alpha_{il}, \beta_i$ and $\delta_i$ are all locally bounded and locally Lipschitz, in the following senses. First, for $i \geq 0$ and $x, y \in \ell_{11}$, we assume that

$$(2.7) \qquad \sum_{l \geq 0, l \neq i} \alpha_{il}(0) \leq a_{00}, \qquad \sum_{l \geq 0, l \neq i}(l+1)\alpha_{il}(0) \leq (i+1)a_{10},$$

$$(2.8) \qquad \sum_{l \geq 0, l \neq i} |\alpha_{il}(x) - \alpha_{il}(y)| \leq a_{01}(x,y)\|x-y\|_1,$$

$$(2.9) \sum_{l \geq 0, l \neq i}(l+1)|\alpha_{il}(x) - \alpha_{il}(y)| \leq (i+1)a_{11}(x,y)\|x-y\|_{11},$$

where the $a_{r0}$ are finite,

$$a_{r1}(x,y) \leq \tilde{a}_{r1}(\|x\|_{11} \wedge \|y\|_{11}), \qquad r = 0, 1,$$

and the $\tilde{a}_{r1}$ are bounded on bounded intervals. Then we assume that, for all $x, y \in \ell_{11}$,

$$\sum_{i\geq 0}(i+1)\beta_i(0) \leq b_{10}, \tag{2.10}$$

$$\sum_{i\geq 0}|\beta_i(x) - \beta_i(y)| \leq b_{01}(x,y)\|x-y\|_1, \tag{2.11}$$

$$\sum_{i\geq 0}(i+1)|\beta_i(x) - \beta_i(y)| \leq b_{11}(x,y)\|x-y\|_{11}, \tag{2.12}$$

where the $b_{r0}$ are finite,

$$b_{r1}(x,y) \leq \tilde{b}_{r1}(\|x\|_{11} \wedge \|y\|_{11}), \qquad r = 0, 1,$$

and the $\tilde{b}_{r1}$ are bounded on bounded intervals; and finally that

$$\sup_{i\geq 0} \delta_i(0) \leq d_0, \tag{2.13}$$

$$\sup_{i\geq 0}|\delta_i(x) - \delta_i(y)| \leq d_1(x,y)\|x-y\|_1, \tag{2.14}$$

where $d_0$ is finite and

$$d_1(x,y) \leq \tilde{d}_1(\|x\|_{11} \wedge \|y\|_{11}),$$

with $\tilde{d}_1$ bounded on finite intervals.

The various assumptions can be understood in the biological context. First, the two norms $\|\cdot\|_1$ and $\|\cdot\|_{11}$ have natural interpretations. The quantity $\|\xi - \eta\|_1$ is the sum of the differences $|\xi^i - \eta^i|$ between the numbers of hosts in states $i = 0, 1, 2, \ldots$ in two host populations $\xi$ and $\eta$; this can be thought of as the natural measure of difference as seen from the hosts' point of view. The corresponding "parasite norm" is then $\|\xi - \eta\|_{11}$, which weights each difference $|\xi^i - \eta^i|$ by the factor $(i+1)$, the number of parasites plus one; in a similar way, writing $x = N^{-1}\xi$, one can interpret $\|x\|_{11}$ as a measure of "parasite density."

The simplest conditions are (2.7) and (2.13), which, together with conditions (2.8) and (2.14) with $y = 0$, ensure that the per capita rates of events involving infection and death are all finite, and bounded by constant multiples of the *host* density $\|x\|_1 + 1$. This is frequently biologically reasonable. For instance, a grazing animal, however hungry, can only consume a limited number of mouthfuls per unit time, an infection event occurring when a mouthful contains infective stages of the parasite (of which there may be many, subject to the remaining conditions). Conditions (2.10) and (2.11) with $y = 0$ can then be interpreted as limiting the overall immigration rate of hosts and the per capita host birth rate. Analogously, conditions (2.8), (2.11) and (2.14) for general $y$ imply that cumulative differences in the above rates between population infection states $x$ and $y$ are limited by multiples of the host norm $\|x-y\|_1$ of the difference between $x$ and $y$, and also that these multiples remain bounded provided that the *smaller* of $\|x\|_{11}$ and $\|y\|_{11}$ remains bounded. In Kretzschmar's (1993) model, and in others in which the per capita infection rate grows as a constant $K$ times the parasite density $\|x\|_{11}$, these conditions are violated, and, without them, the coupling argument of Section 5 fails; however, realistic modifications of these models are covered; see Section 6.

The remaining conditions concern parasite weighted analogues of the preceding conditions. Conditions (2.10) and (2.12) with $y=0$ constrain the overall rate of flow of parasites into the system through immigration to be finite, and bounded if the parasite density remains bounded; condition (2.12) also limits the way in which this influx may depend on the infection state. Conditions (2.7) and (2.9) impose analogous restrictions on the rates at which individual parasites can cause further infection, limiting the output per parasite. Bounds on the parasite weighted quantities are needed to establish the accuracy of approximation that we prove; see the remark following Lemma 4.3. Our choice of bounds allows considerable freedom.

Note that the above discussion relates only to the state-dependent elements in the transition rates. The conditions on the state-independent rates (2.3) and (2.4) have similar effects, but are less restrictive. For instance, the state-independent element in the death rates may increase rather rapidly with parasite burden (2.4), whereas the state-dependent elements are more strongly restricted (2.13), (2.14). The particular form of the conditions on the state-dependent rates is in part dictated by the coupling argument of Section 5, and may well not be the weakest possible for results such as ours to be true. Similarly, the condition (2.5) is of a purely technical nature, though presumably always satisfied in practice.

REMARK. The assumptions made about the $\alpha_{ij}(x)$ and $\beta_i(x)$ have certain general consequences. One is that the total number of hosts has to be finite almost surely for all $t$. This can be seen by comparison with a pure birth process, since the number of hosts $\|X\|_1$ only increases through immigration, and the total rate of immigration $N\sum_{i\geq 0}\beta_i(N^{-1}X)$ does not exceed $Nb_{10} + \tilde{b}_{01}(0)\|X\|_1$. Hence, for any $T > 0$,

$$\mathbf{E}\|X(T)\|_1 \leq N(1 + b_{10}/\tilde{b}_{01})\exp\{T\tilde{b}_{01}\} \tag{2.15}$$

and

$$\lim_{M\to\infty}\mathbf{P}\left[\sup_{0\leq t\leq T}\|X_N(t)\|_1 > NM\right] = 0. \tag{2.16}$$

Now, if $N^{-1}\|X\|_1 \leq M$, it follows from (2.8) that

$$\sum_{l\geq 0}\alpha_{il}(N^{-1}X) \leq a_{00} + \tilde{a}_{01}(0)M,$$

for all $i \geq 0$. Hence, and because $W$ is nonexplosive, it follows that, on the event $\{\sup_{0 \leq t \leq T} \|X_N(t)\|_1 \leq NM\}$, the $X$-chain makes a.s. only finitely many transitions in $[0,T]$. Letting $M \to \infty$, it follows from (2.16) that a.s. only finitely many transitions can occur in the $X$-chain in any finite time interval.

**3. The differential equations.** We assume deterministic initial conditions $X_N(0)$ for each $N$; in fact, because of the Markov property, we could equally well let $X_N(0)$ be measurable with respect to events up to time zero. Our aim is to approximate the evolution of the process $N^{-1}X_N(t)$ when $N$ is large. A natural candidate approximation is given by the solution $x_N$ to the "average drift" infinite dimensional differential equation

$$
\begin{aligned}
\frac{dx^i(t)}{dt} &= \sum_{l \geq 0} x^l(t)\bar{\alpha}_{li} + \sum_{l \neq i} x^l(t)\alpha_{li}(x(t)) - x^i(t)\sum_{l \neq i} \alpha_{il}(x(t)) \\
&\quad + \beta_i(x(t)) - x^i(t)\delta_i(x(t)), \qquad i \geq 0,
\end{aligned}
\tag{3.1}
$$

with initial condition $x_N(0) = N^{-1}X_N(0)$. The following theorem, the main result of the paper, shows that $x_N$ indeed provides a suitable approximation, and quantifies the $\ell_1$-error in the approximation. For convenience, we consider initial conditions that are close to a fixed element $x_0 \in \ell_{11}$.

THEOREM 3.1. *Suppose that (2.3)–(2.14) hold, and that $x_N(0) := N^{-1}X_N(0)$ satisfies $\|x_N(0) - x_0\|_{11} \to 0$ as $N \to \infty$ for some $x_0 \in \ell_{11}$. Let $[0, t_{\max})$ denote the interval on which the equation (3.1) with $x_0$ as initial condition has a solution $x$ belonging to $\ell_{11}$. Then, for any $T < t_{\max}$, there exists a constant $K(T)$ such that, as $N \to \infty$,*

$$\mathbf{P}\bigg[N^{-1}\sup_{0 \leq t \leq T}\|X_N(t) - Nx_N(t)\|_1 > K(T)N^{-1/2}\log^{3/2} N\bigg] = O(N^{-1/2}),$$

*where $x_N$ solves (3.1) with $x_N(0) = N^{-1}X_N(0)$.*

Proving this theorem is the substance of this and the next two sections.

First, it is by no means obvious that equation (3.1) has a solution, something that is only proved in Theorem 4.2. And even if a solution exists, it is not perhaps immediate that it has to belong to the nonnegative cone, in contrast to the case of the stochastic model. To temporarily accommodate this possibility, we extend the definitions of $\alpha_{il}$, $\beta_i$ and $\delta_i$, setting

$$\alpha_{il}(u) = \alpha_{il}(u_+), \qquad \beta_i(u) = \beta_i(u_+), \qquad \delta_i(u) = \delta_i(u_+), \qquad u \in \mathbb{R}_+^\infty,$$

where $u_+^i := \max(u^i, 0)$, $i \geq 0$, and observing that conditions (2.8)–(2.14) are still satisfied, with $x$ and $y$ replaced by $x_+$ and $y_+$, respectively, as arguments of $a_{01}, a_{11}, b_{01}, b_{11}$ and $d_1$.

Temporarily suppressing the $N$-dependence, note that equation (3.1), as in Arrigoni (2003), can be compactly expressed in the form

$$\frac{dx}{dt} = Ax + F(x), \qquad x(0) = N^{-1}X_N(0), \tag{3.2}$$

where $A$ is a linear operator given by

$$(Ax)^i = \sum_{l \geq 0} x^l \bar{\alpha}_{li}, \qquad i \geq 0, \tag{3.3}$$

and $F$ is an operator given by

$$F(x)^i = \sum_{l \neq i} x^l \alpha_{li}(x) - x^i \sum_{l \neq i} \alpha_{il}(x) + \beta_i(x) - x^i \delta_i(x), \qquad i \geq 0. \tag{3.4}$$

If $A$ generates a $C_0$ (strongly continuous) semigroup $T(\cdot)$ on a Banach space $S$, then every solution $x$ of (3.2) with $x(0) \in S$ also satisfies the integral equation

$$x(t) = T(t)x(0) + \int_0^t T(t-s)F(x(s))\,ds. \tag{3.5}$$

A continuous solution $x$ of the integral equation (3.5) is called a *mild* solution of the initial value problem (3.2), and if $F$ is locally Lipschitz continuous and $x(0) \in S$, then (3.5) has a unique continuous solution on some nonempty interval $[0, t_{\max})$ [Pazy (1983), Theorem 1.4, Chapter 6]. In our case, $A^T$ is the infinitesimal matrix of the Markov process $W$, and by standard Markov theory [Kendall and Reuter (1957), pages 111 and 114–115], we can identify $T(t)x$ with $P(t)^T x$ for any $x \in \ell_1$, $T$ being strongly continuous on $\ell_1$. However, in order to establish the bounds that we have claimed, we need to take $S$ to be the space $\ell_{11}$. We therefore show that $T$ is also strongly continuous in $\ell_{11}$, and that $F$ is locally $\ell_{11}$-Lipschitz continuous, from which the existence and uniqueness of a continuous solution in $\ell_{11}$ to the integral equation (3.5) then follows.

LEMMA 3.2.　$T$ *is strongly continuous in* $\ell_{11}$.

PROOF. To prove this lemma, note that every sequence $x$ such that $\|x\|_{11} < \infty$ can be approximated in $\ell_{11}$ by sequences with bounded support, which are all in the domain of $A$, so $D(A)$ is dense in $\ell_{11}$. We now need to check that $T(t)\,\ell_{11} \subseteq \ell_{11}$, and that $T$ is $\ell_{11}$-strongly continuous at 0. First, for every $x$ with $\|x\|_{11} < \infty$, we have $\|P(t)^T x\|_{11} < \infty$ for all times $t$ from (2.3), since

$$\|P(t)^T x\|_{11} = \sum_{j \geq 0}(j+1)\left|\sum_{i \geq 0} x^i p_{ij}(t)\right|$$

(3.6)
$$\leq \sum_{i \geq 0} |x^i|(i+1)e^{wt} = e^{wt}\|x\|_{11} < \infty.$$

For strong continuity, taking any $x \in \ell_{11}$, we have

$$\|P(t)^T x - x\|_{11}$$
$$\leq \sum_{j \geq 0}(j+1)\left\{\sum_{i \neq j} |x^i| p_{ij}(t) + |x^j|(1 - p_{jj}(t))\right\}$$
$$= \left\{\sum_{i \geq 0} |x^i|\{\mathbf{E}_i^0(W(t)+1) - (i+1)\mathbf{P}_i[W(t)=i]\}\right.$$
$$\left. + \sum_{j \geq 0}(j+1)|x^j|(1 - p_{jj}(t))\right\}.$$

Now $\lim_{t \downarrow 0} p_{jj}(t) = 1$ for all $j$, and, by (2.3),

$$0 \leq \mathbf{E}_i^0(W(t)+1) - (i+1)\mathbf{P}_i[W(t)=i] \leq (i+1)e^{wt}$$

and

$$\limsup_{t \downarrow 0} \mathbf{E}_i^0(W(t)+1) \leq \lim_{t \downarrow 0}(i+1)e^{wt} = i+1;$$

$$\liminf_{t \downarrow 0} \mathbf{E}_i^0(W(t)+1) \geq \liminf_{t \downarrow 0} \sum_{j=0}^{i} p_{ij}(t)(j+1) = i+1.$$

Hence, since $\|x\|_{11} < \infty$, $\lim_{t \downarrow 0} \|P(t)^T x - x\|_{11} = 0$ by dominated convergence. □

LEMMA 3.3. *The function $F$ defined in (3.4) is locally Lipschitz continuous in the $\ell_{11}$-norm.*

PROOF. For $x, y \geq 0$ such that $\|x\|_{11}, \|y\|_{11} \leq M$, using assumptions (2.7)–(2.14), we have

$$\|F(x) - F(y)\|_{11}$$
$$\leq \sum_{i \geq 0}(i+1)\sum_{l \neq i} |x^l \alpha_{li}(x) - y^l \alpha_{li}(y)|$$

$$+ \sum_{i\geq 0}(i+1)\left|x^i\sum_{l\neq i}\alpha_{il}(x) - y^i\sum_{l\neq i}\alpha_{il}(y)\right|$$

$$+ \sum_{i\geq 0}(i+1)|\beta_i(x) - \beta_i(y)| + \sum_{i\geq 0}(i+1)|x^i\delta_i(x) - y^i\delta_i(y)|$$

$$\leq \sum_{l\geq 0}|x^l - y^l|\sum_{i\neq l}(i+1)\alpha_{li}(x) + \sum_{l\geq 0}y^l\sum_{i\neq l}(i+1)|\alpha_{li}(x) - \alpha_{li}(y)|$$

$$+ \sum_{i\geq 0}(i+1)|x^i - y^i|\sum_{l\neq i}\alpha_{il}(x) + \sum_{i\geq 0}(i+1)y^i\sum_{l\neq i}|\alpha_{il}(x) - \alpha_{il}(y)|$$

$$+ b_{11}(x,y)\|x - y\|_{11} + \sum_{i\geq 0}(i+1)|x^i - y^i|\delta_i(x)$$

$$+ \sum_{i\geq 0}(i+1)y^i|\delta_i(x) - \delta_i(y)|$$

$$\leq \{a_{10} + \tilde{a}_{11}(0)\|x\|_{11}\}\|x - y\|_{11} + \tilde{a}_{11}(M)\|x - y\|_{11}\|y\|_{11}$$

$$+ \{a_{00} + \tilde{a}_{01}(0)\}\|x\|_1\|x - y\|_{11} + \tilde{a}_{01}(M)\|y\|_{11}\|x - y\|_1$$

$$+ \tilde{b}_{11}(M)\|x - y\|_{11}$$

$$+ \{d_0 + \tilde{d}_1(0)\}\|x\|_1\|x - y\|_{11} + \tilde{d}_1(M)\|x - y\|_1\|y\|_{11}$$

$$\leq F_M\|x - y\|_{11},$$

where

$$\begin{aligned}
(3.7) \quad F_M &:= a_{10} + \tilde{a}_{11}(0)M + M\tilde{a}_{11}(M) \\
&\quad + a_{00} + \tilde{a}_{01}(0)M + M\tilde{a}_{01}(M) \\
&\quad + \tilde{b}_{11}(M) + d_0 + \tilde{d}_1(0)M + M\tilde{d}_1(M). \quad \square
\end{aligned}$$

From these two lemmas, it follows that the differential equation system (3.1) has a unique *weak* solution, so that we at least have a function $x_N$ to give substance to the statement of Theorem 3.1. It is later shown in Theorem 4.2 that, under our conditions, $x_N$ is in fact a classical solution to the system (3.1).

Our main aim is to approximate a single random process $N^{-1}X_N$ by the solution $x_N$ of (3.5) with initial condition $x_N(0) = N^{-1}X_N(0)$, and to give an error bound that is, in principle, computable. However, in order to fix the definitions of the constants appearing in the error bounds, we have framed Theorem 3.1 in terms of a sequence of processes indexed by $N$, assuming that $N^{-1}X_N(0) \to x_0$ in $\ell_{11}$ as $N \to \infty$, for some $x_0 \in \ell_{11}$. It is then natural to be able to approximate all of the processes $N^{-1}X_N$ by the *single* solution $x$ to (3.5) which has $x(0) = x_0$. The next lemma shows that this poses no problems, because the solution of (3.5) depends in locally Lipschitz fashion on the initial conditions, within its interval of existence.

LEMMA 3.4. *Fix a solution $x$ to the integral equation (3.5) with $x(0) \in \ell_{11}$, and suppose that $T < t_{\max}$. Then there is an $\varepsilon > 0$ such that, if $y$ is a solution with initial condition $y(0)$ satisfying $\|y(0) - x(0)\|_{11} \leq \varepsilon$, then*

$$\sup_{0 \leq t \leq T} \|x(t) - y(t)\|_{11} \leq \|x(0) - y(0)\|_{11} C_T,$$

*for a constant $C_T < \infty$.*

PROOF. From the integral equation (3.5) together with (3.6), it follows that, if $\|x(0) - y(0)\|_{11} \leq \varepsilon$, then

$$\|x(t) - y(t)\|_{11} \leq \varepsilon e^{wt} + \int_0^t F_{2M_T} \|x(s) - y(s)\|_{11} e^{w(t-s)} \, ds,$$

where $F_M$ is defined in (3.7) and

$$M_T := \sup_{0 \leq t \leq T} \|x(t)\|_{11},$$

provided also that $\sup_{0 \leq t \leq T} \|y(t)\|_{11} \leq 2M_T$. By Gronwall's inequality, it then follows that

$$\sup_{0 \leq t \leq T} \|x(t) - y(t)\|_{11} \leq \|x(0) - y(0)\|_{11} C_T \leq \varepsilon C_T,$$

for a constant $C_T < \infty$. This implies that $\sup_{0 \leq t \leq T} \|y(t)\|_{11} \leq 2M_T$ is indeed satisfied if $\varepsilon < M_T/C_T$, and the lemma follows. □

COROLLARY 3.5. *Under the conditions of Theorem 3.1, if also $\|x_N(0) - x_0\|_{11} = O(N^{-1/2})$, then $x_N$ can be replaced by $x$ in the statement, without altering the order of the approximation.*

PROOF. Combine Theorem 3.1 with Lemma 3.4. □

It also follows from Lemma 3.4 that, if $\|x_N(0) - x_0\|_{11} \to 0$, then for all $N$ large enough

(3.8) $$\left| \sup_{0 \leq t \leq T} \|x_N(t)\|_{11} - M_T \right| \leq \|x_N(0) - x_0\|_{11} C_T,$$

provided that $T < t_{\max}$. In particular, if $t_{\max}^N$ denotes the maximum time such that $x_N$ is uniquely defined on $[0, t_{\max}^N)$, then $\liminf_{N \to \infty} t_{\max}^N \geq t_{\max}$.

**4. The independent sum approximation.** The next step in proving Theorem 3.1 is to consider an approximating model $\widetilde{X}_N(\cdot)$, starting with $\widetilde{X}_N(0) = X_N(0)$, and consisting of independent individuals. Each individual's parasite load evolves according to a time inhomogeneous Markov process $\widetilde{W}$ on $\mathbb{Z}_+ \cup \Delta$ with infinitesimal matrix defined by

$$
\begin{aligned}
q_{lj}(t) &= \bar{\alpha}_{lj} + \tilde{\alpha}_{lj}(t), & j \neq l, \Delta,\ l \geq 0, \\
q_{ll}(t) &= -\sum_{j \neq l} q_{lj}(t) - \bar{\delta}_l - \tilde{\delta}_l(t), & l \geq 0, \\
q_{l\Delta}(t) &= \bar{\delta}_l + \tilde{\delta}_l(t), & l \geq 0,
\end{aligned}
\tag{4.1}
$$

where

$$
\tilde{\alpha}_{il}(t) := \alpha_{il}(x_N(t)); \qquad \tilde{\delta}_i(t) := \delta_i(x_N(t));
\tag{4.2}
$$

and, for $i, j \in \mathbb{Z}_+ \cup \Delta$, we shall write

$$
\tilde{p}_{ij}(s,t) := \mathbf{P}[\widetilde{W}(t) = j \mid \widetilde{W}(s) = i], \qquad s < t.
\tag{4.3}
$$

In addition, individuals may immigrate, with rates $N\tilde{\beta}_i(t)$, where

$$
\tilde{\beta}_i(t) := \beta_i(x_N(t)).
\tag{4.4}
$$

The process $\widetilde{X}_N$ differs from $X_N$ in having the nonlinear elements of the transition rates made linear, by replacing the Lipschitz state-dependent elements $\alpha_{ij}(x)$, $\beta_i(x)$ and $\delta_i(x)$ at any time $t$ by their "typical" values $\tilde{\alpha}_{ij}(t)$, $\tilde{\beta}_i(t)$ and $\tilde{\delta}_i(t)$. Our strategy will be first to show that the process $\widetilde{X}_N$ stays close to the deterministic process $Nx_N(t)$ with high probability, and then to show that, if this is the case, then $X_N$ also stays close to $\widetilde{X}_N$, again with high probability. However, we shall first use the process $\widetilde{X}_N$ to improve our knowledge about the weak solution $x$ to (3.2).

Let us start by introducing some further notation. Fixing $T < t_{\max}$, define

$$
M_T := M_T(x) := \sup_{0 \leq t \leq T} \sum_{i \geq 1} (i+1)|x^i(t)|;
\tag{4.5}
$$

$$
G_T := G_T(x) := \sup_{0 \leq t \leq T} \sum_{i \geq 0} |x^i(t)|,
\tag{4.6}
$$

and write $M_T^N := M_T(x_N)$, $G_T^N := G_T(x_N)$. Note that $M_T$ is finite if $x(0) \in \ell_{11}$, because the mild solution $x$ is $\ell_{11}$-continuous, and that $M_T^N \geq 1$ and $G_T^N \geq 1$ whenever $\|x_N(0)\|_1 = 1$, as is always the case here.

It is immediate from Lemma 3.4 that if $\|x_N(0) - x(0)\|_{11} \to 0$, then $M_T^N \leq M_T + 1$ for all $N$ large enough. Furthermore, it then follows that

$$
G_T^N \leq M_T^N \leq M_T + 1
\tag{4.7}
$$

for all $N$ sufficiently large. Hence, using also dominated convergence, we deduce that $G_T^N \le G_T + 1$ for all $N$ sufficiently large.

Our first result of the section controls the mean of the process $N^{-1}\widetilde{X}_N$ in the $\ell_{11}$-norm.

LEMMA 4.1. *Under conditions (2.3)–(2.14), for any $X_N(0) \in \ell_{11}$ and any $T < t_{\max}^N$, we have*

$$\sup_{0 \le t \le T} N^{-1} \sum_{l \ge 0} (l+1) \mathbf{E} \widetilde{X}_N^l(t)$$

$$\le \{N^{-1}\|X_N(0)\|_{11} + T(b_{10} + \tilde{b}_{11}(0)M_T^N)\} e^{(w + a_0^* + a_1^* M_T^N)T} < \infty,$$

*where $w$ is as in assumption (2.3), $a_0^* := a_{10} - a_{00}$ and $a_1^* := \tilde{a}_{11}(0) - \tilde{a}_{01}(0)$.*

PROOF. Neglecting the individuals in the cemetery state $\Delta$, the process $\widetilde{X}_N$ can be represented by setting

$$(4.8) \qquad \widetilde{X}_N(t) = \sum_{i \ge 0} \sum_{j=1}^{X_N^i(0)} e(\widetilde{W}_{ij}(t)) + \sum_{i \ge 0} \sum_{j=1}^{R_i(t)} e(\widetilde{W}'_{ij}(t - \tau_{ij})),$$

where $\widetilde{W}_{ij}$ and $\widetilde{W}'_{ij}$, $i \ge 0, j \ge 1$, are independent copies of $\widetilde{W}$, with $\widetilde{W}_{ij}$ and $\widetilde{W}'_{ij}$ starting at $i$, and the $\tau_{ij}$, $j \ge 1$, are the successive event times of independent time inhomogeneous Poisson (counting) processes $R_i$ with rates $N\tilde{\beta}_i(t)$, which are also independent of all the $\widetilde{W}_{ij}$ and $\widetilde{W}'_{ij}$; as usual, $e(l)$ denotes the $l$th coordinate vector. Hence, it follows that, given $X_N(0)$,

$$\sum_{l \ge 0} (l+1) \mathbf{E} \widetilde{X}_N^l(t)$$

$$= \sum_{i \ge 0} \left\{ X_N^i(0) \mathbf{E}_i^0 \{\widetilde{W}(t) + 1\} + N \int_0^t \tilde{\beta}_i(u) \mathbf{E}_i^0 \{\widetilde{W}(t-u) + 1\} \, du \right\},$$

where $\mathbf{E}_i^0$ is as defined for (2.3).

Now $\widetilde{W}$ has paths which are piecewise paths of $W$, but with extra killing because of the $\tilde{\delta}_i(u)$ components of the rates, and with extra jumps, $\alpha$-jumps, say, because of the $\tilde{\alpha}_{ij}(u)$ components. The killing we can neglect, since it serves only to reduce $\mathbf{E}_i^0 \widetilde{W}(t)$. For the remainder, by assumption (2.8), the rate of occurrence is at most $\chi := \{a_{00} + \tilde{a}_{01}(0)M_T^N\}$, irrespective of state and time. So, defining

$$c_{ij}(u) := \tilde{\alpha}_{ij}(u)/\chi, \qquad i,j \ge 0, \ j \ne i;$$

$$c_{ii}(u) := 1 - \sum_{j \ne i} c_{ij}(u), \qquad i \ge 0,$$

we can construct the $\alpha$-jumps by taking them to occur at the event times of a Poisson process $R$ of rate $\chi$, with jump distribution for a jump at time $u$ given by $c_i.(u)$ if $\widetilde{W}(u-) = i$. Note that, in this case, no jump is realized with probability $c_{ii}(u)$. Conditional on $R$ having events at times $0 < t_1 < \cdots < t_r < t$ between $0$ and $t$, we thus have

$$\mathbf{E}_i^0\{\widetilde{W}(t) + 1 \mid t_1, \ldots, t_r\}$$
$$\leq \sum_{j_1 \geq 0} p_{ij_1}(t_1) \sum_{l_1 \geq 0} c_{j_1 l_1}(t_1) \sum_{j_2 \geq 0} p_{l_1 j_2}(t_2 - t_1) \sum_{l_2 \geq 0} c_{j_2 l_2}(t_2) \cdots$$
$$\cdots \sum_{j_r \geq 0} p_{l_{r-1} j_r}(t_r - t_{r-1}) \sum_{l_r \geq 0} c_{j_r l_r}(t_r) \sum_{j \geq 0} (j+1) p_{l_r j}(t - t_r),$$

where, as before, $p_{ij}(t) = \mathbf{P}[W(t) = j \mid W(0) = i]$. Applying assumptions (2.3), (2.7) and (2.9) to the last two sums, we have

$$\sum_{l_r \geq 0} c_{j_r l_r}(t_r) \sum_{j \geq 0} (j+1) p_{l_r j}(t - t_r) \leq \sum_{l_r \geq 0} c_{j_r l_r}(t_r)(l_r + 1) e^{w(t - t_r)}$$
$$\leq a_3^N (j_r + 1) e^{w(t - t_r)},$$

where $a_3^N := \{a_{10} + \tilde{a}_{11}(0) M_T^N\}/\chi$. It thus follows that

$$\mathbf{E}_i^0\{\widetilde{W}(t) + 1 \mid t_1, \ldots, t_r\} \leq a_3^N \mathbf{E}_i^0\{\widetilde{W}(t_r) + 1 \mid t_1, \ldots, t_{r-1}\} e^{w(t - t_r)}.$$

Arguing inductively, this implies that

$$\mathbf{E}_i^0\{\widetilde{W}(t) + 1 \mid t_1, \ldots, t_r\} \leq (i+1)\{a_3^N\}^r e^{wt},$$

and hence, unconditionally, that

$$\mathbf{E}_i^0\{\widetilde{W}(t) + 1\} \leq (i+1) e^{wt} \mathbf{E}\{(a_3^N)^{R(t)}\}$$
(4.9)
$$\leq (i+1) \exp\{(w + (a_3^N - 1)\chi)t\}.$$

The remainder of the proof is immediate. □

Armed with this estimate, we can now proceed to identify $N^{-1} \mathbf{E} \widetilde{X}_N(t)$ with $x_N(t)$, at the same time proving that the mild solution $x_N$ is in fact a classical solution to the infinite differential equation (3.1) with initial condition $N^{-1} X_N(0)$.

First, define the "linearized" version of (3.1):

$$\frac{dy^i(t)}{dt} = \sum_{l \geq 1} y^l(t) \bar{\alpha}_{li} + \sum_{l \geq 0} y^l(t) \tilde{\alpha}_{li}(t) - y^i(t) \sum_{l \geq 0} \tilde{\alpha}_{il}(t)$$
(4.10)
$$+ \tilde{\beta}_i(t) - y^i(t) \tilde{\delta}_i(t), \qquad i \geq 0,$$

where $\tilde{\alpha}$, $\tilde{\beta}$ and $\tilde{\delta}$ are as in (4.2) and (4.4), to be solved in $t \in [0, T]$ for an unknown function $y$. Clearly these equations have $x_N$ itself as a mild solution in $\ell_{11}$, and, by Pazy (1983), Theorem 3.2, Chapter 6, the mild solution is unique under our assumptions, since now

$$\widetilde{F}(t, u)^i := \sum_{l \geq 0} u^l \tilde{\alpha}_{li}(t) - u^i \sum_{l \geq 0} \tilde{\alpha}_{il}(t) + \tilde{\beta}_i(t) - u^i \tilde{\delta}_i(t)$$

is $\ell_{11}$ locally Lipschitz in $u \in S$ with constant

$$a_{00} + a_{10} + d_0 + M_T^N \{\tilde{a}_{01}(0) + \tilde{a}_{11}(0) + \tilde{d}_1(0)\}$$

and

$$\|\widetilde{F}(s, u) - \widetilde{F}(t, u)\|_{11}$$
$$\leq (M_T^N \{\tilde{a}_{11}(M_T^N) + \tilde{a}_{01}(M_T^N) + \tilde{d}_1(M_T^N)\} + \tilde{b}_{11}(M_T^N))\|x_N(s) - x_N(t)\|_{11},$$

whenever $\|u\|_{11} \leq M_T^N$, so that $\widetilde{F}$ is $t$-uniformly continuous on bounded intervals (contained in $[0, t_{\max}^N)$), because $x_N$ is. We now show that $y(t) = N^{-1}\mathbf{E}\widetilde{X}_N(t)$ solves (4.10), and indeed as a *classical* solution. Since it also therefore solves (3.5), and since this equation has a unique solution, it follows that $y$ is the same as $x_N$, and that it is the classical solution to equation (3.1) with initial condition $N^{-1}X_N(0)$.

THEOREM 4.2. *Under conditions (2.3)–(2.14), for any fixed $X_N(0) \in \ell_{11}$, the function $y(t) := N^{-1}\mathbf{E}\widetilde{X}_N(t)$ satisfies the system (4.10) with initial condition $N^{-1}X_N(0)$ on any interval $[0, T]$ with $T < t_{\max}^N$. It is hence the unique classical solution $x_N$ to (3.1) for this initial condition.*

PROOF. Let $\widetilde{X}_{N1}^j(t)$ denote the number of particles present at time 0 that are still present and in state $j$ at time $t$; let $\widetilde{X}_{N2}^j(t)$ denote the number of particles that immigrated after time 0 and are present and in state $j$ at time $t$. Then

$$\begin{aligned}(4.11) \quad \mathbf{E}\widetilde{X}_N^j(t) &= \mathbf{E}\widetilde{X}_{N1}^j(t) + \mathbf{E}\widetilde{X}_{N2}^j(t) \\ &= \sum_{i \geq 0} X_N^i(0)\tilde{p}_{ij}(0, t) + N\int_0^t \sum_{i \geq 0} \tilde{\beta}_i(u)\tilde{p}_{ij}(u, t)\,du,\end{aligned}$$

where $\tilde{p}_{ij}(u, v)$ is as defined in (4.3). Note that the expectations are finite, since, from conditions (2.10) and (2.11) and by (4.7), uniformly in $u \in [0, T]$,

$$(4.12) \qquad \sum_{i \geq 0} \tilde{\beta}_i(u) \leq b_{10} + \tilde{b}_{01}(0) G_T^N < \infty.$$

Then, defining $q_{jl}(t)$ as in (4.1), it follows for each $j \geq 0$ that the quantities $-q_{jj}(t)$ are bounded, uniformly in $t$, since by conditions (2.13), (2.14) and (2.8),

$$(4.13) \qquad \tilde{\delta}_j(t) \leq d_0 + \tilde{d}_1(0)G_T^N; \qquad \sum_{l \geq 0} \tilde{\alpha}_{jl}(t) \leq a_{00} + \tilde{a}_{01}(0)G_T^N,$$

and since $\alpha^*(j) = \sum_{l \neq j} \bar{\alpha}_{jl}$ is finite. Note further that, from the forward equations,

$$\frac{\partial}{\partial t}\tilde{p}_{ij}(u,t) = \sum_{l \geq 0} \tilde{p}_{il}(u,t)q_{lj}(t);$$

see Iosifescu and Tautu (1973), Corollary to Theorem 2.3.8, page 214.

Now we have

$$\mathbf{E}\widetilde{X}_{N1}^j(t) = \sum_{i \geq 0} X_N^i(0)\tilde{p}_{ij}(0,t)$$

$$(4.14) \qquad = X_N^j(0) + \sum_{i \geq 0} X_N^i(0) \int_0^t \sum_{l \geq 0} \tilde{p}_{il}(0,u)q_{lj}(u)\,du$$

$$= X_N^j(0) + \int_0^t \sum_{l \geq 0} \mathbf{E}\widetilde{X}_{N1}^l(u)q_{lj}(u)\,du,$$

with no problems about reordering, because of the uniform boundedness discussed above, and since only one of the $q_{lj}$ is negative. Then also, with rearrangements similarly justified, we define

$$Q_t := N\int_0^t \left\{ \int_0^v \sum_{i \geq 0} \tilde{\beta}_i(u) \sum_{l \geq 0} \tilde{p}_{il}(u,v)q_{lj}(v)\,du \right\} dv;$$

taking the $i$-sum first, and then the $u$-integral, we obtain

$$Q_t = \int_0^t \sum_{l \geq 0} \mathbf{E}\widetilde{X}_{N2}^l(v)q_{lj}(v)\,dv;$$

taking the $l$-sum first, we have

$$Q_t = N\int_0^t \left\{ \int_0^t \sum_{i \geq 0} \tilde{\beta}_i(u)\frac{\partial}{\partial v}\tilde{p}_{ij}(u,v)\mathbf{1}_{[0,v]}(u)\,du \right\} dv$$

$$= N\int_0^t \sum_{i \geq 0} \tilde{\beta}_i(u)\{\tilde{p}_{ij}(u,t) - \tilde{p}_{ij}(u,u)\}\,du$$

$$= \mathbf{E}\widetilde{X}_{N2}^j(t) - N\int_0^t \tilde{\beta}_j(u)\,du.$$

From these two representations of $Q_t$, it follows that

$$(4.15) \quad \mathbf{E}\widetilde{X}_{N2}^j(t) = \int_0^t \left\{ N\tilde{\beta}_j(u) + \sum_{l\geq 0} \mathbf{E}\widetilde{X}_{N2}^l(u) q_{lj}(u) \right\} du;$$

combining (4.15) and (4.15), we thus have

$$(4.16) \quad N^{-1}\mathbf{E}\widetilde{X}_N^j(t)$$
$$= N^{-1}X_N^j(0) + \int_0^t \left\{ \tilde{\beta}_j(u) + \sum_{l\geq 0} N^{-1}\mathbf{E}\widetilde{X}_N^l(u) q_{lj}(u) \right\} du.$$

Since the right-hand side is an indefinite integral up to $t$, it follows that $N^{-1}\mathbf{E}\widetilde{X}_N^j(t)$ is continuous in $t$, for each $j$. The quantities $q_{jl}(t)$ are all continuous, because $x_N$ is $\ell_{11}$-continuous in $t$ and the $\alpha_{il}(x)$ and $\delta_l(x)$ are $\ell_{11}$-Lipschitz, and also, for $q_{ll}(t)$, from assumption (2.8). Then, from Lemma 4.1, we also have

$$\sum_{j\geq J} \mathbf{E}\widetilde{X}_N^j(t) \leq (J+1)^{-1}\{\|X_N(0)\|_{11} + NT(b_{10} + \tilde{b}_{11}(0)M_T^N)\}e^{(w+a_0^* + a_1^* M_T^N)T},$$

so that, in view of assumption (2.5) and of (4.13), the sum on the right-hand side of (4.16) is uniformly convergent, and its sum continuous. Hence (4.16) can be differentiated with respect to $t$ to recover the system (4.10), proving the theorem. □

Our next result shows that, under appropriate conditions, $N^{-1}\widetilde{X}_N(t)$ and $x(t)$ are close in $\ell_1$-norm at any *fixed* $t$, with very high probability.

LEMMA 4.3. *Suppose that conditions (2.3)–(2.14) are satisfied, that $X_N(0) \in \ell_{11}$ and that $N \geq 9$. Then, for any $t \in [0,T]$ with $T < t_{\max}^N$,*

$$\mathbf{E}\|\widetilde{X}_N(t) - Nx_N(t)\|_1 \leq 3(M_T^N + 1)\sqrt{N\log N}.$$

*Furthermore, for any $r > 0$, there exist constants $K_r^{(1)} > 1$ and $K_r^{(2)}$ such that*

$$\mathbf{P}[\|\widetilde{X}_N(t) - Nx_N(t)\|_1 > K_r^{(1)}(M_T^N + 1)N^{1/2}\log^{3/2} N] \leq K_r^{(2)} G_T^N N^{-r}.$$

PROOF. For a sum $W$ of independent indicator random variables with mean $M$, and for any $\delta > 0$, it follows from the Chernoff inequalities that

$$(4.17) \quad \max\{\mathbf{P}[W > M(1+\delta)], \mathbf{P}[W < M(1-\delta)]\} \leq \exp\{-M\delta^2/(2+\delta)\};$$

see Chung and Lu (2006), Theorem 4. Now the quantity $\widetilde{X}_N^j(t)$ can be expressed as a sum of independent random variables $Y_1^j, \ldots, Y_N^j$ and $Y'$, where $Y_k^j$ is the indicator of the event that the $k$th initial individual is in state $j$ at time $t$, and $Y'$ is a Poisson random variable with mean $N \int_0^t \sum_{i \geq 0} \beta_i(s) \tilde{p}_{ij}(s,t) \, ds$. Hence, by the simple observation that $\mathrm{Var}(\widetilde{X}_N^i(t)) \leq \mathbf{E}(\widetilde{X}_N^i(t)) = Nx_N^i(t)$, we have

$$\mathbf{E}|\widetilde{X}_N^i(t) - Nx_N^i(t)| \leq \sqrt{Nx_N^i(t)} \wedge \{2Nx_N^i(t)\}, \tag{4.18}$$

and, by (4.17), for any $a \geq 2$ and $N \geq 3$, we have

$$\mathbf{P}[|\widetilde{X}_N^i(t) - Nx_N^i(t)| > a\sqrt{Nx_N^i(t) \log N}] \leq 2N^{-a/2},$$

so long as $Nx_N^i(t) \geq 1$. By Lemma A.1(iii), in view of (4.18), it now follows immediately that

$$\mathbf{E}\|\widetilde{X}_N(t) - Nx_N(t)\|_1 \leq 3(M_T^N + 1)\sqrt{N \log N},$$

giving the first statement.

For the second, let $I_N(t) := \{i : x_N^i(t) \geq 1/N\}$. Then it is immediate that $|I_N(t)| \leq NG_T^N$, so that, if

$$B_N(t) := \bigcap_{i \in I_N(t)} \{|\widetilde{X}_N^i(t) - Nx_N^i(t)| \leq a\sqrt{Nx_N^i(t) \log N}\},$$

then

$$\mathbf{P}[B_N^C(t)] \leq 2G_T^N N^{1-a/2}. \tag{4.19}$$

On the other hand, on the event $B_N(t)$, it follows from Lemma A.1 (i) that

$$\sum_{i \in I_N(t)} |\widetilde{X}_N^i(t) - Nx_N^i(t)| \leq a \log N \sum_{i \in I(t)} \sqrt{Nx_N^i(t)}$$
$$\leq a \log N \, (M_T^N + 1)\sqrt{N \log N}. \tag{4.20}$$

For the remaining indices, we note that $S_N(t) := \sum_{i \notin I_N(t)} \widetilde{X}_N^i(t)$ is also a sum of many independent indicator random variables plus an independent Poisson component. Using (4.17), we thus have

$$\mathbf{P}\left[S_N(t) > \sum_{i \notin I_N(t)} Nx_N^i(t) + N^{1/2}(M_T^N + 1)\sqrt{\log N}\right]$$
$$\leq \exp\{-N^{1/2}(M_T^N + 1)\sqrt{\log N}/3\} \leq \exp\{-N^{1/2}/3\}, \tag{4.21}$$

since $\delta := N^{1/2}(M_T^N + 1)\sqrt{\log N} / \sum_{i \notin I_N(t)} Nx_N^i(t) \geq 1$ from Lemma A.1(ii); otherwise, we have

$$\sum_{i \notin I_N(t)} \widetilde{X}_N^i(t) \leq 2N^{1/2}(M_T^N + 1)\sqrt{\log N}, \tag{4.22}$$

again from Lemma A.1 (ii). Now, fixing any $r > 0$ and taking $a = 2(r+1)$, the second part of the lemma follows from (4.19)–(4.22). □

REMARK. The argument above makes essential use of the finiteness of $M_T^N$, through Lemma A.1. For $M_T^N < \infty$, we needed the solution $x_N$ to (3.5) to be $\ell_{11}$-continuous, and hence needed to consider equation (3.5) with respect to $\ell_{11}$. Now, for the accuracy of approximation given in Lemma 4.3, the condition $M_T^N < \infty$ is in fact close to being necessary. To see this, note that, if $Z_i \sim \text{Po}(Np_i)$ for $i \geq 1$, for some choice of $p_i$, then

$$\sum_{i:p_i \geq 1/N} \mathbf{E}|Z_i - Np_i| \asymp N^{1/2} \sum_{i:p_i \geq 1/N} \sqrt{p_i}.$$

In the arguments above, we have the $x_N^i(t)$ for the $p_i$, with the $\widetilde{X}_N^i(t)$, which are close to being Poisson distributed, in place of $Z_i$, and we use the fact that $\sum_{i \geq 1} ip_i < \infty$. Suppose instead that $p_i \sim ci^{-1-\eta}$ for some $c > 0$ and $0 < \eta < 1$, so that $\sum_{i \geq 1} i^\eta p_i = \infty$. Then it follows that

$$\sum_{i:p_i \geq 1/N} \sqrt{p_i} \asymp N^{(1-\eta)/2(1+\eta)},$$

and hence, that

$$\sum_{i \geq 1} \mathbf{E}|Z_i - Np_i| \geq KN^{1/2}N^{(1-\eta)/2(1+\eta)}.$$

Thus, for such $p_i$, an approximation as close as that of Lemma 4.3 cannot be attained, because the mean $\ell_1$-distance would be at least of order as big as $N^\gamma$ for $\gamma = 1/(1+\eta) > 1/2$. Note that such circumstances would arise in our model, if we took, for instance, $\bar{\alpha}_{i,i-1} = i$, $i \geq 2$, $\bar{\delta}_1 = 1$, $\beta_i(x) = i^{-1-\eta}$ for all $x$, and set all other elements of the transition rates equal to zero. The resulting stochastic model has no interactions, so that the processes $\widetilde{X}_N$ and $X_N$ have the same distribution, that of a particular multitype Markov immigration–death process. At stationarity, the mean number $p_i$ of individuals of type $i$ satisfies $p_i \sim ci^{-1-\eta}$ as $i \to \infty$, for $c = 1/\eta$. Of course, with this choice, the $\beta_i(0)$ violate condition (2.10).

The next lemma is used to control the fluctuations of $\widetilde{X}_N$ between close time points. We define the quantity

$$H_T^N := 2^{m_2-1}m_1 + \{b_{10} + a_{00} + \tilde{b}_{01}(0)$$
$$+ d_0 + G_T^N(\tilde{a}_{01}(0) + \tilde{d}_1(0))\}/\lceil NM_T^N \rceil^{m_2-1},$$

which will be used as part of an upper bound for the transition rates of the process $\widetilde{X}_N$ on $[0,T]$, noting that

(4.23) $$1 \leq H_T^N \leq H^* M_T^N,$$

where $H^*$ does not depend on $T$ or $N$.

LEMMA 4.4. *Suppose that conditions (2.3)–(2.14) are satisfied. Then, if $h \le 1/(2\lceil NM_T^N \rceil^{m_2} H_T^N)$, $t \le T - h$, and if $\|\widetilde{X}_N(t) - Nx_N(t)\|_1 \le KN^{1/2}\log^{3/2} N$, it follows that*

$$\mathbf{P}\left[\sup_{0 \le u \le h} \|\widetilde{X}_N(t+u) - \widetilde{X}_N(t)\|_1 > KN^{1/2}\log^{3/2} N + a\log N\right] \le N^{-a/6},$$

*for any $a \ge 2$ and $N \ge 3$.*

PROOF. At time $t$, there are $\|\widetilde{X}_N(t)\|_1$ individuals in the system, each of which evolves independently of the others over the interval $s \in [t, t+h]$; in addition, new immigrants may arrive. During the interval $[t, t+h]$, an individual in state $i \ge 0$ at time $t$ has probability

$$\exp\left\{-\left(h(\alpha^*(i) + \bar{\delta}_i) + \int_0^h \delta_i(x(t+u))\,du + \int_0^h \sum_{l \ne i}\alpha_{il}(x(t+u))\,du\right)\right\}$$

of not changing state; and the expected number of immigrants is Poisson distributed with mean

$$N\sum_{i \ge 0}\int_0^h \beta_i(x(t+u))\,du.$$

Now consider the total number $R(t,h)$ of individuals that either change state or immigrate during the interval $[t, t+h]$. For each $i \ge 0$, the individuals in state $i$ at time $t$ can be split into two groups, the first containing $\widetilde{X}_N^i(t) \wedge Nx_N^i(t)$ randomly chosen individuals, and the second containing the remainder. Adding over $i$, the numbers in the second group add up to at most $KN^{1/2}\log^{3/2} N$, by assumption. Then, from the observations above, the mean number of individuals in the first group that change state in $[t, t+h]$ is at most

$$\sum_{i=0}^{2\lceil NM_T^N \rceil - 1} hNx_N^i(t)m_1(i+1)^{m_2} + \sum_{i \ge 2\lceil NM_T^N \rceil} Nx_N^i(t)$$

(4.24)
$$+ \sum_{i \ge 0} hNx_N^i(t)\{d_0 + \tilde{d}_1(0)G_T^N + a_{00} + \tilde{a}_{01}(0)G_T^N\}$$

from (2.4) and (2.7)–(2.14). Finally, the expected number of immigrants in $[t, t+h]$ is at most

(4.25) $$hN\{b_{10} + \tilde{b}_{01}(0)G_T^N\}.$$

Adding (4.24) and (4.25), and recalling the assumption on $h$, we obtain an expected number of events in these categories of at most

$$hNM_T^N m_1\{2\lceil NM_T^N \rceil\}^{m_2-1} + \tfrac{1}{2}$$

$$(4.26) \quad + hN\{b_{10} + G_T^N(d_0 + a_{00} + \tilde{b}_{01}(0)) + [G_T^N]^2(\tilde{d}_1(0) + \tilde{a}_{01}(0))\}$$
$$\leq \tfrac{1}{2} + h\lceil NM_T^N\rceil^{m_2} H_T^N \leq 1.$$

Applying the Chernoff bounds (4.17), the probability that more than $a \log N \geq 2$ events of these kinds occur in $[t, t+h]$ is thus at most $N^{-a/6}$, implying in sum that

$$\mathbf{P}[R(t,h) \geq KN^{1/2} \log^{3/2} N + a\log N] \leq N^{-a/6}.$$

Since $\sup_{0 \leq u \leq h} \|\widetilde{X}_N(t+u) - \widetilde{X}_N(t)\|_1 \leq R(t,h)$, the lemma follows. $\square$

We are now in a position to prove the main result of the section, showing that the independent sum process $N^{-1}\widetilde{X}_N$ is a good approximation to $x_N$, uniformly in $[0,T]$.

THEOREM 4.5. *Under conditions (2.3)–(2.14), and for $T < t_{\max}^N$, there exist constants $1 \leq K_r^{(3)}, K_r^{(4)} < \infty$ for each $r > 0$ such that, for $N$ large enough,*

$$\mathbf{P}\left[\sup_{0 \leq t \leq T} \|N^{-1}\widetilde{X}_N(t) - x_N(t)\|_1 > K_r^{(4)}(M_T^N + F_{M_T^N})N^{-1/2}\log^{3/2} N\right]$$
$$\leq K_r^{(3)}(M_T^N)^{m_2+2} N^{-r},$$

*where $F_M$ is as in (3.7).*

PROOF. Suppose that $N \geq 9$. Divide the interval $[0,T]$ into $\lceil 2T\lceil NM_T^N\rceil^{m_2} \times H_T^N\rceil$ intervals $[t_l, t_{l+1}]$ of lengths $h_l = t_{l+1} - t_l \leq 1/(2\lceil NM_T^N\rceil^{m_2} H_T^N)$. Apply Lemma 4.3 with $r + m_2$ for $r$ and with $t = t_l$ for each $l$, and apply Lemma 4.4 with $a = 6(r + m_2)$ and with $t = t_l$ and $h = h_l$ for each $l$; except on a set of probability at most

$$\lceil 2T\lceil NM_T^N\rceil^{m_2} H_T^N\rceil(K_{r+m_2}^{(2)} G_T^N N^{-r-m_2} + N^{-r-m_2}),$$

we have

$$\sup_{0 \leq t \leq T} \|N^{-1}\widetilde{X}_N(t) - x_N(t)\|_1$$
$$(4.27) \quad \leq 2K_{r+m_2}^{(1)}(M_T^N + 1)N^{-1/2}\log^{3/2} N + 6(r+m_2)N^{-1}\log N$$
$$+ \sup_{0 \leq s,t \leq T; |s-t| \leq 1/(2\lceil NM_T^N\rceil^{m_2} H_T^N)} \|x_N(s) - x_N(t)\|_1.$$

Now, since $x_N$ satisfies (3.5), it follows that, for $0 \leq u \leq h$,

$$\|x_N(t+u) - x_N(t)\|_1 \leq \|x_N(t)P(u) - x_N(t)\|_1$$
$$+ \int_0^u \|F(x_N(t+v))P(h-v)\|_1 \, dv.$$

By (3.6), we have

$$\|F(x_N(t+v))P(h-v)\|_{11} \le e^{wu}\|F(x_N(t+v))\|_{11} \le e^{wh}M_T^N F_{M_T^N},$$

the last inequality following from Lemma 3.3, so that therefore

$$\int_0^u \|F(x_N(t+v))P(h-v)\|_1 \, dv \le he^{wh}M_T^N F_{M_T^N} \le \frac{e^{wT}F_{M_T^N}}{2N}$$

for $h \le 1/\{2\lceil NM_T^N \rceil^{m_2} H_T^N\}$, since $m_2 \ge 1$ and $H_T^N \ge 1$. Then

$$\|x_N(t)P(u) - x_N(t)\|_1 \le 2\sum_{j\ge 0} |x_N^j(t)|(1 - p_{jj}(u)),$$

and, by part of the calculation in (4.26), we find that

$$\sum_{j\ge 0} |x_N^j(t)|(1 - p_{jj}(u)) \le u \sum_{j=0}^{2\lceil NM_T^N \rceil - 1} |x_N^j(t)|m_1(j+1)^{m_2} + \sum_{j\ge 2\lceil NM_T^N \rceil} |x_N^j(t)|$$

$$\le hm_1(2\lceil NM_T^N \rceil^{m_2-1})M_T^N + \frac{1}{2N},$$

which is at most $3/(2N)$ if $h \le 1/(2\lceil NM_T^N \rceil^{m_2} H_T^N)$.

Hence, the third term in (4.27) is of order $(1 + F_{M_T^N})N^{-1}$ under the conditions of the theorem, and the result follows. $\square$

**5. The main approximation.** We now turn to estimating the deviations of the process $\widetilde{X}_N$ from the actual process $X_N$ of interest. We do so by coupling the processes in such a way that the "distance" between them cannot increase too much over any finite time interval. In our coupling, we pair each individual in state $i \ge 1$ in $X_N(0)$ with a corresponding individual in state $i$ in $\widetilde{X}_N(0)$ so that all their $\mu$- and $\bar{\delta}$-transitions are identical. This process entails an implicit labeling, which we suppress from the notation. Now the remaining transitions have rates which are not quite the same in the two processes, and hence, the two can gradually drift apart. Our strategy is to make their transitions identical as far as we can, but, once a transition in one process is not matched in the other, the individuals are decoupled thereafter. For our purposes, it is simply enough to show that the *number* of decoupled pairs is sufficiently small; what pairs of states these individuals occupy is immaterial.

We realize the coupling between $X_N$ and $\widetilde{X}_N$ in terms of a four component process $Z(\cdot)$ with

$$Z(t) = ((Z_l^i(t), i \ge 0, 1 \le l \le 3), Z_4(t)) \in \mathcal{X}^3 \times \mathbb{Z}_+,$$

constructed in such a way that we can define $X_N(\cdot) = Z_1(\cdot) + Z_2(\cdot)$ and $\widetilde{X}_N(\cdot) = Z_1(\cdot) + Z_3(\cdot)$, and starting with $Z_1(0) = X_N(0) = \widetilde{X}_N(0)$, $Z_2(0) = Z_3(0) = 0 \in \mathcal{X}$, and $Z_4(0) = 0$. The component $Z_4$ is used only to keep count of certain uncoupled individuals: either unmatched $Z_2$-immigrants, or $Z_3$ individuals that die, or $Z_2$ individuals created at the death of (one member of) a coupled pair. The transition rates of $Z$ are given as follows, using the notation $e_l(i)$ for the $i$th coordinate vector in the $l$th copy of $\mathcal{X}$, and writing $X = Z_1 + Z_2$. For the $\bar\alpha$- and $\alpha$-transitions, at time $t$ and for any $i \neq l$, we have

$$Z \to Z + (e_1(l) - e_1(i)) \text{ at rate } Z_1^i\{\bar\alpha_{il} + (\alpha_{il}(N^{-1}X) \wedge \alpha_{il}(x_N(t)))\};$$
$$Z \to Z + (e_2(l) + e_3(i) - e_1(i)) \text{ at rate } Z_1^i\{\alpha_{il}(N^{-1}X) - \alpha_{il}(x_N(t))\}^+;$$
$$Z \to Z + (e_2(i) + e_3(l) - e_1(i)) \text{ at rate } Z_1^i\{\alpha_{il}(N^{-1}X) - \alpha_{il}(x_N(t))\}^-;$$
$$Z \to Z + (e_2(l) - e_2(i)) \text{ at rate } Z_2^i\{\bar\alpha_{il} + \alpha_{il}(N^{-1}X)\};$$
$$Z \to Z + (e_3(l) - e_3(i)) \text{ at rate } Z_3^i\{\bar\alpha_{il} + \alpha_{il}(x_N(t))\},$$

with possibilities for individuals in the two processes to become uncoupled, when $N^{-1}X \neq x(t)$. For the immigration transitions, we have

$$Z \to Z + e_1(i) \text{ at rate } N\{\beta_i(N^{-1}X) \wedge \beta_i(x_N(t))\}, \qquad i \geq 0;$$
$$Z \to Z + e_2(i) + e_4 \text{ at rate } N\{\beta_i(N^{-1}X) - \beta_i(x_N(t))\}^+, \qquad i \geq 0;$$
$$Z \to Z + e_3(i) \text{ at rate } N\{\beta_i(N^{-1}X) - \beta_i(x_N(t))\}^-, \qquad i \geq 0,$$

with some immigrations not being precisely matched; the second transition includes an $e_4$ to ensure that each individual in $Z_2$ has a counterpart in either $Z_3$ or $Z_4$. For the deaths, we have

$$Z \to Z - e_1(i) \text{ at rate } Z_1^i\{\bar\delta_i + (\delta_i(N^{-1}X) \wedge \delta_i(x_N(t)))\}, \qquad i \geq 0;$$
$$Z \to Z - e_1(i) + e_3(i) \text{ at rate } Z_1^i\{\delta_i(N^{-1}X) - \delta_i(x_N(t))\}^+, \qquad i \geq 0;$$
$$Z \to Z - e_1(i) + e_2(i) + e_4 \text{ at rate } Z_1^i\{\delta_i(N^{-1}X) - \delta_i(x_N(t))\}^-, \qquad i \geq 0;$$
$$Z \to Z - e_2(i) \text{ at rate } Z_2^i\{\bar\delta_i + \delta_i(N^{-1}X)\}, \qquad i \geq 0;$$
$$Z \to Z - e_3(i) + e_4 \text{ at rate } Z_3^i\{\bar\delta_i + \delta_i(x_N(t))\}, \qquad i \geq 0,$$

where $Z_4(\cdot)$ is also used to count the deaths of uncoupled $Z_3$-individuals, and uncoupled deaths in $\tilde X_N$ of coupled $Z_1$ individuals. With this joint construction, we have arranged that

$$(5.1) \qquad \sum_{i \geq 0} Z_2^i(t) \leq Z_4(t) + \sum_{i \geq 0} Z_3^i(t)$$

for all $t$, and that

$$(5.2) \qquad V_N(t) := Z_4(t) + \sum_{i \geq 0} Z_3^i(t)$$

is a counting process. We allow unmatched deaths in the $Z_2$-process. We thus have the bound

$$\|X_N(t) - \widetilde{X}_N(t)\|_1 = \|(Z_1(t) + Z_2(t)) - (Z_1(t) + Z_3(t))\|_1$$
$$\leq \sum_{i \geq 0} \{Z_2^i(t) + Z_3^i(t)\} \quad (5.3)$$
$$\leq 2\left\{Z_4(t) + \sum_{i \geq 0} Z_3(t)\right\} = 2V_N(t),$$

for all $t$, by (5.1).

Now $V_N$ has a compensator $A_N$ with intensity $a_N$, satisfying

$$a_N(t) = \sum_{i \geq 0} Z_1^i(t) \sum_{l \geq 0} |\alpha_{il}(N^{-1}X_N(t)) - \alpha_{il}(x_N(t))|$$
$$+ N \sum_{i \geq 0} |\beta_i(N^{-1}X_N(t)) - \beta_i(x_N(t))|$$
$$+ \sum_{i \geq 0} Z_1^i(t)|\delta_i(N^{-1}X_N(t)) - \delta_i(x_N(t))|$$
$$\leq \sum_{i \geq 0} \widetilde{X}_N^i(t) \sum_{l \geq 0} |\alpha_{il}(N^{-1}X_N(t)) - \alpha_{il}(x_N(t))|$$
$$+ N \sum_{i \geq 0} |\beta_i(N^{-1}X_N(t)) - \beta_i(x_N(t))|$$
$$+ \sum_{i \geq 0} \widetilde{X}_N^i(t)|\delta_i(N^{-1}X_N(t)) - \delta_i(x_N(t))|.$$

Now, condition (2.8) implies that, uniformly in $i$,

$$\sum_{l \geq 0} |\alpha_{il}(N^{-1}X_N(t)) - \alpha_{il}(x_N(t))| \leq \tilde{a}_{01}(\|x_N(t)\|_{11})\|N^{-1}X_N(t) - x_N(t)\|_1.$$

Hence,

$$N^{-1}a_N(t) \leq \left(\sum_{i \geq 0} x_N^i(t)\tilde{a}_{01}(M_T^N) + \tilde{b}_{01}(M_T^N) + \sum_{i \geq 0} x_N^i(t)\tilde{d}_1(M_T^N)\right)$$
$$\times \|N^{-1}X_N(t) - x_N(t)\|_1$$
$$+ \|N^{-1}\widetilde{X}_N(t) - x_N(t)\|_1(\tilde{a}_{01}(M_T^N) + \tilde{d}_1(M_T^N))$$
$$\times \|N^{-1}X_N(t) - x_N(t)\|_1$$
$$\leq \{H_T^{(1,N)} + H_T^{(2,N)}\|N^{-1}\widetilde{X}_N(t) - x_N(t)\|_1\}\|N^{-1}X_N(t) - x_N(t)\|_1,$$

where
$$H_T^{(1,N)} = G_T^N \tilde{a}_{01}(M_T^N) + \tilde{b}_{01}(M_T^N) + G_T^N \tilde{d}_1(M_T^N);$$
$$H_T^{(2,N)} = \tilde{a}_{01}(M_T^N) + \tilde{d}_1(M_T^N).$$

In particular, defining
$$\tau_N := \inf\{t \geq 0 : \|N^{-1}\widetilde{X}_N(t) - x_N(t)\|_1 \geq 1\},$$

it follows that

(5.4)
$$N^{-1}a_N(t \wedge \tau_N)$$
$$\leq \{H_T^{(1,N)} + H_T^{(2,N)}\}\|N^{-1}X_N(t \wedge \tau_N) - x_N(t \wedge \tau_N)\|_1.$$

For the next two lemmas, we shall restrict the range of $N$ in a way that is asymptotically unimportant. We shall suppose that $N$ satisfies the inequalities

(5.5)
$$K_{r_0}^{(4)}(M_T^N + F_{M_T^N})N^{-1/2}\log^{3/2} N \leq 1;$$
$$N > \max\{(K_{r_0}^{(3)})^{1/(m_2+2)} M_T^N, 9\},$$

where $r_0 = m_2 + 2$, and the quantities $K_r^{(3)}$ and $K_r^{(4)}$ are as for Theorem 4.5.

LEMMA 5.1. *Under conditions (2.3)–(2.14), for any $t \in [0, T]$ and for all $N$ satisfying (5.5), we have*

$$N^{-1}\mathbf{E}\|X_N(t \wedge \tau_N) - \widetilde{X}_N(t \wedge \tau_N)\|_1$$
$$\leq 14 M_T^N N^{-1/2}\sqrt{\log N}\, t(H_T^{(1,N)} + H_T^{(2,N)})\exp\{2t(H_T^{(1,N)} + H_T^{(2,N)})\}.$$

PROOF. Write $\mathcal{M}_N(\cdot) := V_N(\cdot) - A_N(\cdot)$. Then, because also

(5.6)
$$\|N^{-1}X_N(t) - x_N(t)\|_1$$
$$\leq \|N^{-1}X_N(t) - N^{-1}\widetilde{X}_N(t)\|_1 + \|N^{-1}\widetilde{X}_N(t) - x_N(t)\|_1,$$

and using (5.4), we have

$$(2N)^{-1}\|X_N(t \wedge \tau_N) - \widetilde{X}_N(t \wedge \tau_N)\|_1$$
$$\leq N^{-1}V(t \wedge \tau_N)$$
$$\leq N^{-1}\mathcal{M}_N(t \wedge \tau_N)$$
(5.7)
$$+ \int_0^{t \wedge \tau_N} \{H_T^{(1,N)} + H_T^{(2,N)}\}\{N^{-1}\|X_N(s \wedge \tau_N) - \widetilde{X}_N(s \wedge \tau_N)\|_1$$
$$+ \|N^{-1}\widetilde{X}_N(s \wedge \tau_N) - x_N(s \wedge \tau_N)\|_1\}\,ds$$

$$\leq N^{-1}\mathcal{M}_N(t\wedge\tau_N)$$
$$+\int_0^t\{H_T^{(1,N)}+H_T^{(2,N)}\}\{N^{-1}\|X_N(s\wedge\tau_N)-\widetilde{X}_N(s\wedge\tau_N)\|_1$$
$$+\|N^{-1}\widetilde{X}_N(s\wedge\tau_N)-x_N(s\wedge\tau_N)\|_1\}\,ds.$$

Now $\mathcal{M}_N(\cdot\wedge\tau_N)$ is a martingale, since, by (5.4),

$$\mathbf{E}\Big\{\int_0^{t\wedge\tau_N}a_N(s)\,ds\Big\}$$
$$\leq\{H_T^{(1,N)}+H_T^{(2,N)}\}\mathbf{E}\Big\{\int_0^{t\wedge\tau_N}\|X_N(s)-Nx_N(s)\|_1\,ds\Big\}$$
$$\leq\{H_T^{(1,N)}+H_T^{(2,N)}\}\int_0^t\{\mathbf{E}\|X_N(s)\|_1+N\|x_N(s)\|_1\}\,ds\ <\ \infty,$$

this last using (2.15). Taking expectations, it thus follows that

$$(2N)^{-1}\mathbf{E}\|X_N(t\wedge\tau_N)-\widetilde{X}_N(t\wedge\tau_N)\|_1$$
$$(5.8)\quad \leq\int_0^t\{H_T^{(1,N)}+H_T^{(2,N)}\}\{N^{-1}\mathbf{E}\|X_N(s\wedge\tau_N)-\widetilde{X}_N(s\wedge\tau_N)\|_1$$
$$+\mathbf{E}\|N^{-1}\widetilde{X}_N(s\wedge\tau_N)-x_N(s\wedge\tau_N)\|_1\}\,ds.$$

Now we have

$$\mathbf{E}\|N^{-1}\widetilde{X}_N(s\wedge\tau_N)-x_N(s\wedge\tau_N)\|_1$$
$$=\mathbf{E}\{\|N^{-1}\widetilde{X}_N(s)-x_N(s)\|_1 I[\tau_N\geq s]\}$$
$$(5.9)\quad +\mathbf{E}\{\|N^{-1}\widetilde{X}_N(\tau_N)-x_N(\tau_N)\|_1 I[\tau_N<s]\}$$
$$\leq\mathbf{E}\|N^{-1}\widetilde{X}_N(s)-x_N(s)\|_1+(1+1/N)\mathbf{P}[\tau_N<s]$$
$$\leq\mathbf{E}\|N^{-1}\widetilde{X}_N(s)-x_N(s)\|_1+(1+1/N)\mathbf{P}[\tau_N<T].$$

For $N$ satisfying (5.5), the first term in (5.9) is bounded by $3(M_T^N+1)N^{-1/2}\times\sqrt{\log N}$ by Lemma 4.3. Also, for such $N$, the event $\{\tau_N<T\}$ lies in the exceptional set for Theorem 4.5 with $r=r_0$, implying that

$$(5.10)\qquad \mathbf{P}[\tau_N<T]\leq K_{r_0}^{(3)}(M_T^N)^{m_2+2}N^{-r_0},$$

so that the second term is no larger than $M_T^N N^{-1/2}\sqrt{\log N}$ if

$$N\ >\ \max\{(K_{r_0}^{(3)}(M_T^N)^{m_2+1})^{1/(r_0-1/2)},9\},$$

which is also true if (5.5) is satisfied. This implies that, for such $N$,

$$(5.11)\quad \mathbf{E}\|N^{-1}\widetilde{X}_N(s\wedge\tau_N)-x_N(s\wedge\tau_N)\|_1\leq 7M_T^N N^{-1/2}\sqrt{\log N}.$$

Using (5.11) in (5.8) and applying Gronwall's inequality, the lemma follows.
$\square$

LEMMA 5.2. *Under conditions (2.3)–(2.14), for any $t \in [0, T]$ and $y > 0$, and for all $N$ satisfying (5.5), we have*

$$\mathbf{P}\left[\sup_{0 \le s \le t} |N^{-1} \mathcal{M}_N(s \wedge \tau_N)| \ge y\right] \le g(t(H_T^{(1,N)} + H_T^{(2,N)})) \, M_T^N y^{-2} N^{-3/2} \sqrt{\log N},$$

*where $g(x) := 7xe^{2x}$.*

PROOF. Since $V_N(\cdot)$ is a counting process with continuous compensator $A_N$, we have from (5.4) and (5.6) that

$$\mathbf{E}\mathcal{M}_N^2(t \wedge \tau_N) = \mathbf{E}A_N(t \wedge \tau_N)$$
$$\le \{H_T^{(1,N)} + H_T^{(2,N)}\} \int_0^t \{\mathbf{E}\|X_N(s \wedge \tau_N) - \widetilde{X}_N(s \wedge \tau_N)\|_1$$
$$+ \mathbf{E}\|\widetilde{X}_N(s \wedge \tau_N) - Nx_N(s \wedge \tau_N)\|_1\} \, ds.$$

The first expectation is bounded using Lemma 5.1, the second from (5.11), from which it follows that

$$\mathbf{E}\mathcal{M}_N^2(t \wedge \tau_N) \le g(t(H_T^{(1,N)} + H_T^{(2,N)})) \, M_T^N N^{1/2} \sqrt{\log N}.$$

The lemma now follows from the Lévy–Kolmogorov inequality. □

We are finally in a position to complete the proof of Theorem 3.1.

PROOF OF THEOREM 3.1. Suppose that $N$ satisfies (5.5). Returning to the inequality (5.7), we can now write

$$(2N)^{-1} \sup_{0 \le s \le t} \|X_N(s \wedge \tau_N) - \widetilde{X}_N(s \wedge \tau_N)\|_1 \le N^{-1} V(t \wedge \tau_N)$$
$$\le N^{-1} \mathcal{M}_N(t \wedge \tau_N)$$
$$+ \int_0^t \{H_T^{(1,N)} + H_T^{(2,N)}\}\{N^{-1}\|X_N(s \wedge \tau_N) - \widetilde{X}_N(s \wedge \tau_N)\|_1$$
$$+ \|N^{-1} \widetilde{X}_N(s \wedge \tau_N) - x_N(s \wedge \tau_N)\|_1\} \, ds.$$

From Lemma 5.2, taking $y = y_N = M_T^N N^{-1/2} \sqrt{\log N}$, we can bound the martingale contribution uniformly on $[0, T]$ by $M_T^N N^{-1/2} \sqrt{\log N}$, except on an event of probability at most

$$g(T(H_T^{(1,N)} + H_T^{(2,N)})) N^{-1/2}.$$

By Theorem 4.5, for any $r > 0$, we can find a constant $K_r$ such that

$$\sup_{0 \le t \le T} \|N^{-1} \widetilde{X}_N(t) - x_N(t)\|_1 \le K_r(M_T^N + F_{M_T^N}) N^{-1/2} \log^{3/2} N,$$

except on an event of probability $O((M_T^N)^{m_2+2}N^{-r})$. Hence, once again by Gronwall's inequality, it follows that, except on these exceptional events,

$$N^{-1}\sup_{0\leq s\leq t}\|X_N(s\wedge\tau_N)-\widetilde{X}_N(s\wedge\tau_N)\|_1$$

$$\leq 2N^{-1/2}e^{2t(H_T^{(1,N)}+H_T^{(2,N)})}$$
$$\times\{M_T^N\sqrt{\log N}+T(H_T^{(1,N)}+H_T^{(2,N)})K_r(M_T^N+F_{M_T^N})\log^{3/2}N\}.$$

Combining this with Theorem 4.5, and since also, by (5.10), $\mathbf{P}[\tau_N<T]=O((M_T^N)^{m_2+2}N^{-r})$ for any $r$, the theorem follows:

$$\mathbf{P}\bigg[N^{-1}\sup_{0\leq t\leq T}\|X_N(t)-Nx_N(t)\|_1 > K(T)N^{-1/2}\log^{3/2}N\bigg]=O(N^{-1/2}).$$

Note that the inequalities (5.5) are satisfied for all $N$ sufficiently large, and that the constant $K(T)$ and the implied constant in $O(N^{-1/2})$ can be chosen uniformly in $N$, because, under the conditions of the theorem, $\|x_N(0)-x(0)\|_{11}\to 0$ as $N\to\infty$, with the result that, for all large enough $N$, $G_T^N$ and $M_T^N$ can be replaced by $G_T+1$ and $M_T+1$ respectively, with the corresponding modifications in $H_T^{(1,N)}$, $H_T^{(2,N)}$ and $F_{M_T^N}$. □

**6. Examples.** In this section we show that the assumptions that we have made are satisfied by a number of epidemic models. These include the models introduced in Barbour and Kafetzaki (1993) and in Barbour (1994), both of which were generalized and studied in depth in Luchsinger (1999, 2001a, 2001b).

6.1. *Luchsinger's nonlinear model.* In Luchsinger's nonlinear model, the total population size is fixed at $N$ at all times, with $\beta_i(x)=\delta_i(x)=\bar{\delta}_i=0$ for all $i\geq 0$ and $x\in\ell_{11}$. The matrix $\bar{\alpha}$ is the superposition of the infinitesimal matrices of a linear pure death process with rate $\mu>0$ and of a catastrophe process which jumps from any state to 0 at constant rate $\kappa\geq 0$. The first of the above expresses the assumption that parasites die independently at rate $\mu$. The second corresponds to the fact that hosts die independently at rate $\kappa$, and their parasites with them; whenever a host dies, it is instantly replaced by a healthy individual. Thus, the positive elements of $\bar{\alpha}$ are given by

$$\bar{\alpha}_{i,i-1}=i\mu,\qquad \bar{\alpha}_{i0}=\kappa,\qquad i\geq 2;\qquad \bar{\alpha}_{10}=\mu+\kappa,$$

and $\bar{\alpha}_{ii}$, $i\geq 1$, is determined by (2.1). The elements $\bar{\alpha}_{0j}$ are all zero. It is easy to check that assumptions (2.3)–(2.5) are satisfied, with $w=0$, $m_1=\mu+\kappa$ and $m_2=1$.

As regards infection, hosts make potentially infectious contacts at rate $\lambda > 0$, and infection can only occur in a currently uninfected host. If a host carrying $i$ parasites contacts a healthy one, infection with $l$ parasites is developed by the healthy host with probability $p_{il}$, where $\sum_{l \geq 0} p_{il} = 1$ for all $i$ and $p_{00} = 1$. Here, the distribution $F_i = (p_{il}, l \geq 0)$ is the $i$-fold convolution of $F_1$, modeling the assumption that, at such a contact, the parasites act independently in transmitting offspring to the previously healthy host. These rules are incorporated by taking

$$\alpha_{0l}(x) = \lambda \sum_{i \geq 1} x^i p_{il}, \qquad l \geq 1, \ x \in \ell_{11},$$

and the remaining $\alpha_{il}(x)$ are all zero. Thus, for $z \geq 0$, we can take $a_{00} = a_{10} = 0$ and

$$\tilde{a}_{01}(z) = \lambda, \qquad \tilde{a}_{11}(z) = \lambda \max\{\theta, 1\},$$

where $\theta$ is the mean of $F_1$, the mean number of offspring transmitted by a parasite during an infectious contact: thus, $\sum_{l \geq 0} p_{il}(l+1) = i\theta + 1$.

6.2. *Luchsinger's linear model.* In Luchsinger's linear model there is tacitly assumed to be an *infinite* pool of potential infectives, so that the 0-coordinate is not required, and its value may if desired be set to 0; the population of interest consists of the infected hosts, whose number may vary. The matrix $\bar{\alpha}$ is the infinitesimal matrix of a simple death process with rate $\mu > 0$, but now restricted to the reduced state space, giving the positive elements

$$\bar{\alpha}_{i,i-1} = i\mu, \qquad i \geq 2;$$

hosts losing infection are now incorporated by using the $\bar{\delta}_i$, with

$$\bar{\delta}_i = \kappa, \qquad i \geq 2; \qquad \bar{\delta}_1 = \kappa + \mu,$$

again with $\bar{\alpha}_{ii}$, $i \geq 2$, determined by (2.1). Assumptions (2.3)–(2.5) are again easily satisfied. Only a member of the pool of uninfected individuals can be infected, and infections with $i$ parasites occur at a rate $\lambda \sum_{l \geq 1} X^l p_{li}$, so that we have

$$\beta_i(x) = \lambda \sum_{l \geq 1} x^l p_{li}, \qquad i \geq 1,$$

with all the $\alpha_{il}(x)$ and $\delta_i(x)$ equal to zero. Here, for $z > 0$, we can take $b_{10} = 0$ and

$$\tilde{b}_{01}(z) = \lambda, \qquad \tilde{b}_{11}(z) = \lambda \max\{\theta, 1\}.$$

6.3. *Kretzschmar's model.* In Kretzschmar's model mortality of parasites is modeled as in Luchsinger's nonlinear model. In addition, hosts can both die and give birth, with death rates increasing with parasite load, and birth rates decreasing. In our formulation, the positive elements of $\bar{\alpha}$ are given by

$$\bar{\alpha}_{i,i-1} = i\mu, \qquad i \geq 1,$$

and we take

$$\bar{\delta}_i = \kappa + i\alpha, \qquad \delta_i(x) = 0,$$

for nonnegative constants $\kappa$ and $\alpha$; again, the $\bar{\alpha}_{ii}$ are determined by (2.1). Assumptions (2.3)–(2.5) are easily satisfied. Hosts are born free of parasites, so that $\beta_i(x) = 0$ for $i \geq 1$, and

$$\beta_0(x) = \beta \sum_{i \geq 0} x^i \xi^i,$$

for some $\beta > 0$ and $0 \leq \xi \leq 1$. Here, we can take $b_{10} = 0$ and, for $z \geq 0$, $b_{01}(z) = b_{11}(z) = \beta$.

Infection only takes place a single parasite at a time, but at a complicated state dependent rate. The $\alpha_{il}(x)$ are all zero except for $l = i + 1$, when

$$\alpha_{i,i+1}(x) = \lambda \varphi(x), \qquad i \geq 0, \qquad \text{with } \varphi(x) = \sum_{j \geq 1} j x^j \Big/ \Big\{ c + \sum_{j \geq 0} x^j \Big\}.$$

Unfortunately, for $c > 0$, we cannot improve substantially on the bound

$$|\varphi(x) - \varphi(y)| \leq c^{-1} \|x - y\|_{11} + c^{-2}(\|x\|_{11} \wedge \|y\|_{11})\|x - y\|_1,$$

which does not yield a suitable bound for $a_{01}(x, y)$, because of the appearance of $\|x - y\|_{11}$. However, the same average rate of parasite transmission is obtained if we define

$$(6.1) \qquad \alpha_{i,i+j}(x) = \lambda x^j \Big/ \Big\{ c + \sum_{l \geq 0} x^l \Big\}, \qquad j \geq 1, \ i \geq 0,$$

or, more generally, much as in Luchsinger's models,

$$(6.2) \qquad \alpha_{i,i+j}(x) = \nu \sum_{l \geq 1} x^l p_{lj} \Big/ \Big\{ c + \sum_{l \geq 0} x^l \Big\}, \qquad j \geq 1, \ i \geq 0,$$

where $\sum_{j \geq 1} j p_{lj} = l\theta$ for each $l \geq 1$, and $\nu\theta = \lambda$. With these rates, for $c > 0$, our conditions are satisfied with

$$a_{00} = a_{10} = 0,$$
$$a_{01}(x, y) = 2\nu c^{-1},$$
$$a_{11}(x, y) = 2\lambda\{c^{-1} + c^{-2}(\|x\|_{11} \wedge \|y\|_{11})\}.$$

In Kretzschmar's model parasites are assumed to be ingested singly, but potentially arbitrarily fast. In our variant there is a fixed rate $\nu$ of taking mouthfuls. Each mouthful leads to the ingestion of $j$ parasites with probability related to the relative frequencies of $l$-hosts, $l \geq 1$, in the population, and to the chance $p_{lj}$ that an $l$-host transmits $j$ parasites in the mouthful. At least for grazing animals, this would seem to be a more plausible description of what might happen. Kretzschmar's model can, for instance, be interpreted as the limit of this model as $\nu \to \infty$, with $\theta = \lambda/\nu$ for fixed $\lambda$ and $\{p_{lj}, j \geq 1\}$ the probabilities from the Poisson distribution $\text{Po}(l\theta)$. Note that, as is to be expected, the value of $a_{01}(x, y)$ tends to infinity in this limit.

## APPENDIX

In this section we prove a lemma used in the proof of Lemma 4.3.

LEMMA A.1. *Suppose that $u_i \geq 0$ for $i \geq 0$, and that $\sum_{i \geq 0}(i+1)u_i \leq M < \infty$ for some $M \geq 1$. Let $I_N := \{i \geq 0 : u_i \geq 1/N\}$. Then, for any $N \geq 9$, we have the following:*

(i) $\quad \displaystyle\sum_{i \in I_N} \sqrt{u_i} \leq (M+1)\sqrt{\log N},$

(ii) $\quad \displaystyle\sum_{i \notin I_N} u_i \leq (M+1)N^{-1/2}\sqrt{\log N},$

(iii) $\quad \displaystyle\sum_{i \geq 0} \{\sqrt{Nu_i} \wedge (2Nu_i)\} \leq 3(M+1)\sqrt{N \log N}.$

PROOF. Define $J_N := \{i \geq 0 : (i+1)\sqrt{u_i \log N} \geq 1\}$ and $K_N := \{i : 0 \leq i < \sqrt{N}\}$. Then we have

$$\sum_{i \in I_N} \sqrt{u_i} = \sum_{i \in I_N \setminus J_N} \sqrt{u_i} + \sum_{i \in I_N \cap J_N} \sqrt{u_i}$$

$$\leq \sum_{(i+1) \leq \sqrt{N/\log N}} \frac{1}{(i+1)\sqrt{\log N}} + \sum_{i \geq 0} (i+1)u_i \sqrt{\log N}$$

$$\leq (M+1)\sqrt{\log N}.$$

A similar calculation then gives

$$\sum_{i \notin I_N} u_i = \sum_{i \in K_N \setminus I_N} u_i + \sum_{i \notin (I_N \cup K_N)} u_i$$

$$\leq \sum_{i=0}^{\lfloor \sqrt{N} \rfloor} N^{-1} + \sum_{i > \lfloor \sqrt{N} \rfloor} N^{-1/2} i u_i$$

$$\leq N^{-1/2} + N^{-1} + MN^{-1/2} \leq N^{-1/2}(M+1)\sqrt{\log N}$$

in $N \geq 9$. Part (iii) combines the two results. $\square$

**Acknowledgments.** A. D. Barbour wishes to thank both the Institute for Mathematical Sciences of the National University of Singapore and Monash University School of Mathematics for providing a welcoming environment while part of this work was accomplished. M. J. Luczak wishes to thank the University of Zürich for their warm hospitality while the work on this project commenced. The authors express their grateful appreciation of the work put into the paper by the referees.

Angewandte Mathematik
Universität Zürich
Winterthurertrasse 190
CH-8057 Zürich
Switzerland
E-mail: a.d.barbour@math.uzh.ch

London School of Economics
Houghton Street
London WC2A 2AE
United Kingdom
E-mail: m.j.luczak@lse.ac.uk